\documentclass[final,3p,times]{elsarticle}
\usepackage{amsfonts}

 \usepackage{graphics}
 \usepackage{graphicx}
 \usepackage{epsfig}

\usepackage{amssymb}
\usepackage{amsmath,latexsym}
\usepackage{mathrsfs}
\usepackage{bm}
\usepackage{amssymb}
\usepackage{hyperref}
\usepackage{hyperref}
\newtheorem{thm}{Theorem}
\newtheorem{lem}[thm]{Lemma}
\newtheorem{pro}[thm]{Proposition}
\newdefinition{de}[thm]{Definition}
\newdefinition{rmk}{Remark}
\newproof{pf}{Proof}
\newproof{pot}{Proof of the existence of Theorem 3.2}
\newtheorem{hypo}{Hypothesis}

\journal{SAIJOPASP}

\begin{document}

\begin{frontmatter}



\title{Backward Stochastic Differential Equations \\
and Feynman-Kac Formula for Multidimensional L\'{e}vy Processes,
with Applications in Finance\tnoteref{t1}} \tnotetext[t1]{This
research was supported by the National Basic Research Program of
China (973 Program) (Program No.2007 CB814903) and the National
Natural Science Foundation of China (Program No.70671069).}

\author[sjtu]{Jianzhong Lin\corref{cor1}}

\ead{jzlin@sjtu.edu.cn} \cortext[cor1]{Corresponding author.}

\address[sjtu]{Department of Mathematics, Shanghai Jiaotong University, Shanghai 200240, China}

\begin{abstract}
In this paper we show the existence and form uniqueness of a solution for
multidimensional backward stochastic differential equations driven
by a multidimensional L\'{e}vy process with moments of all orders.
The results are important from a pure mathematical point of view as
well as in the world of finance: an application to Clark-Ocone and
Feynman-Kac formulas for multidimensional L\'{e}vy processes is
presented. Moreover, the Feynman-Kac formula and the related partial
differential integral equations provide an analogue of the famous
Black-Scholes partial differential equation and thus can be used for
the purpose of option pricing in a multidimensional L\'{e}vy market.
\end{abstract}


\end{frontmatter}


\section{Introduction}\label{sec1}

A linear version of Backward Stochastic Differential Equations
(BSDEs in short) driven by Brownian motion was initially consided by
Bismut(1973)in the context of optimal control. Nonlinear BSDEs were
later introduced by Pardoux and Peng (1990) and independently by
Duffie and Epstein(1992). The BSDE theory has found wide
applications in partial differential equation theory, stochastic
controls and, particularly, mathematical finance(see El Karoui and
Quenez(1997); Ma and Yong(1999)).

Situ (1997) studied BSDEs driven by a Brownian motion and a Poisson
point process. Ouknine (1998) considered BSDEs driven by a Poisson
random measure. Nualart and Schoutens (2000) proved a martingale
representation theorem for L\'{e}vy processes satisfying some
exponential moment condition. By using this martingale
representation result, Nualart and Schoutens (2001) established the
existence and uniqueness of solutions for BSDEs driven by a L\'{e}vy
process of the kind considered in Nualart and Schoutens (2000),
furthermore Nualart and Schoutens (2001) presented the Clark-Ocone
and the Feynman-Kac formulas, the related Partial Differential
Integral Equation (PDIE) and their applications in finance.

In the past twenty years, there is already a growing interest for
multidimensional L\'{e}vy Processes. Some concepts and basic
properties about multidimensional L\'{e}vy Processes were summarized
in Sato (1999). Applications of multidimensional L\'{e}vy Processes
to analyzing biomolecular (DNA and protein) data and one-server
light traffic queues were explored by Dembo, Karlin and Zeittouni
(1994). A small deviations property of multidimensional L\'{e}vy
Processes were discussed by Simon (2003). In finance research,
practically all financial applications require a multidimensional
model with dependence between components: examples are basket option
pricing, portfolio optimization, simulation of risk scenarios for
portfolios. In most of these applications, jumps in the price
process must be taken into account. Cont and Tankov (2004)
systematically investigated these problems in multidimensional
L\'{e}vy market. In addition, the optimal portfolios in
multidimensional L\'{e}vy market is discussed by Emmer and
Kl\"{u}ppelberg (2004), and option pricing is investigated by
Reich,N., Schwab, C. and Winter, C.(2009). Some simulation
approaches for multidimensional L\'{e}vy processes are also
investigated in Cohen and Rosi\'{n}ski (2007). L\'{e}vy copulas was
also suggested by Kallsen and Tankov (2006) in order to characterize
the dependence among components of multidimensional L\'{e}vy
Processes.

Recently, a martingale representation theorem for multidimensional
L\'{e}vy processes was also proved in Lin (2011), and the obtained
representation formula was similar as that in Nualart and Schoutens
(2000). The purpose of this paper is to use this martingale
representation result obtained in Lin (2011) to establish the
existence and form uniqueness of solutions for BSDE's driven by a
multivariate L\'{e}vy process considered in Lin (2011). Although the
proof techniques are similar to those in Nualart and Schoutens
(2001), the results are important from a pure mathematical point of
view as well as in the world of finance. This is illustrated in the
applications. The resulting Clark-Ocone and Feynman-Kac formulas are
fundamental ingredients in the build up of an Malliavin calculus for
multidimensional L\'{e}vy processes. Moreover, the Feynman-Kac
formula and the related Partial Differential Integral Equation
(PDIE) also have an important application in finance: they provide
us an analogue of the famous Black-Scholes partial differential
equation and is used for the purpose of option pricing in a
multidimensional L\'{e}vy market.

The paper is organized as follows. Section 2 contains some
preliminaries on multidimensional L\'{e}vy processes. Section 3
contains the main result on BSDEs driven by multidimensional L'{e}vy
processes. In Section 4 we have included some applications of BSDE's
driven by multidimensional L\'{e}vy processes to the Clark-Ocone,
the Feynman-Kac formulas, and option pricing in a multivariate
L\'{e}vy market. Finally, in the appendix one can find detailed
proofs of the main results.

\section{Preliminary}\label{sec2}

A $\mathbb{R}^n$-valued stochastic process
$X=\{X(t)=(X_1(t),X_2(t),\cdots,X_n(t))',t\geq 0\}$ defined in
complete probability space $(\Omega,\mathscr{F},\mathbb{P})$ is
called \textsl{L\'{e}vy process} if $X$ has stationary and
independent increments and $X(0)=\bm{0}$. A L\'{e}vy process
possesses a c\`{a}dl\`{a}g modification and we will always assume
that we are using this c\`{a}dl\`{a}g version. If we let
$\mathscr{F}_t=\mathscr{G}_t\vee\mathscr{N}$, where
$\mathscr{G}_t=\sigma\{X(s),0\leq s\leq t\}$ is the natural
filtration of $X$, and $\mathscr{N}$ are the $\mathbb{P}-$null sets
of $\mathscr{F}$, then $\{\mathscr{F}_t,t\geq 0\}$ is a right
continuous family of $\sigma-$fields. We assume that $\mathscr{F}$
is generated by $X$. For an up-to-date and comprehensive account of
L\'{e}vy processes we refer the reader to Bertoin (1996) and Sato
(1999).

Let $X$ be a L\'{e}vy process and denote by
\begin{eqnarray}
X(t-)=\lim\limits_{s\rightarrow t,s<t}X(s),\quad t>0 ,\nonumber
\end{eqnarray}
the left limit process and by $\triangle X(t)=X(t)-X(t-)$ the jump
size at time $t$. It is known that the law of $X(t)$ is
\textsl{infinitely divisible} with characteristic function of the
form
\begin{eqnarray}
E\left[exp(i\bm{\theta}\cdot
X(t))\right]=\left(\phi(\bm{\theta})\right)^t,\quad
\bm{\theta}=(\theta_1,\theta_2,\cdots,\theta_n)\in
\mathbb{R}^n\nonumber
\end{eqnarray}
where $\phi(\bm{\theta})$ is the characteristic function of
$\bm{X}(1)$. The function $\psi(\bm{\theta})=log\phi(\bm{\theta})$
is called the \textsl{characteristic exponent} and it satisfies the
following famous L\'{e}vy-Khintchine formula (Bertoin, 1996):
\begin{eqnarray}
\psi(\bm{\theta})=-\frac{1}{2}\bm{\theta}\cdot
\Sigma\bm{\theta}+\textrm{i}\bm{a}\cdot
\bm{\theta}+\int_{\mathbb{R}^n}\left(exp(i\bm{\theta}\cdot
\bm{x})-1-\textrm{i}\bm{\theta}\cdot \bm{x} 1_{|\bm{x}|\leq
1}\right)\nu(d\bm{x}).\nonumber
\end{eqnarray}
where $\bm{a},\bm{x}\in \mathbb{R}^n$, $\Sigma$ is a symmetric
nonnegative-definite $n\times n$ matrix, and $\nu$ is a measure on
$\mathbb{R}^n\backslash\{o\}$ with $\int(\|\bm{x}\|^2\wedge
1)\nu(d\bm{x})<\infty$. The measure $\nu$ is called the
\textsl{L\'{e}vy measure} of $X$.

Throughout this paper, we will use the standard multi-index
notation. We denote by $\mathbb{N}_0$ the set of nonnegative
integers. A multi-index is usually denoted by $\bm{p}$,
$\bm{p}=(p_1,p_2,\cdots,p_n)\in\mathbb{N}_0^n$. Whenever $\bm{p}$
appears with subscript or superscript, it means a multi-index. In
this spirit, for example, for $\bm{x}=(x_1,\cdots,x_n)$, a monomial
in variables $x_1,\cdots,x_n$ is denoted by
$\bm{x}^{\bm{p}}=x_1^{p_1}\cdots x_n^{p_n}$.  In addition, we also
define $\bm{p}!=p_1!\cdots p_n!$ and $|\bm{p}|=p_1+\cdots+p_n$; and
if $\bm{p}$, $\bm{q}\in\mathbb{N}_0^n$, then we define
$\delta_{\bm{p},\bm{q}}=\delta_{\bm{p}_1,\bm{q}_1}\cdots\delta_{\bm{p}_n,\bm{q}_n}$.

\begin{hypo}
We will suppose in the remaining of
the paper that the L\'{e}vy measure satisfies for some
$\varepsilon>0$, and $\lambda>0$,
\begin{eqnarray}
\int_{|\bm{x}|\geq \epsilon}exp(\lambda
\|\bm{x}\|)\nu(d\bm{x})<\infty. \nonumber
\end{eqnarray}
\end{hypo}

This implies that
\begin{eqnarray}
\int \bm{x}^{\bm{p}}\nu(d\bm{x})<\infty. \quad |\bm{p}|\geq
2\nonumber
\end{eqnarray}
and that the characteristic function $E\left[exp(i\bm{\theta}\cdot
X(t))\right]$ is analytic in a neighborhood of origin $\bm{o}$. As a
consequence, $X(t)$ has moments of all orders and the polynomials
are dense in $L^2(\mathbb{R}^n,\mathbb{P}\circ X(t)^{-1})$ for all
$t>0$.

Fix a time interval $[0,T]$ and set
$L_T^2=L^2(\Omega,\mathscr{F}_T,\mathbb{P})$. We will denote by
$\mathscr{P}$ the predictable sub-$\sigma$-field of
$\mathscr{F}_T\otimes\mathscr{B}_{[0,T]}$. First we introduce some
notation:
\begin{itemize}
\item [$\bullet$]: Let $H_T^2$ denote the space of square integrable
and $\mathscr{F}_t-$progressively one-dimensional measurable
processes $\phi=\{\phi(t),t\in[0,T]\}$ such that
\begin{eqnarray}
\|\phi\|^2=\mathbb{E}\left[\int_0^T\|\phi(t)\|^2dt\right]<\infty.\nonumber
\end{eqnarray}
\item [$\bullet$]: $M_T^2$ will denote the subspace of $H_T^2$
formed by predictable processes.
\item [$\bullet$]: $(H_T^2(l^2))^m$ and $(M_T^2(l^2))^m$ are the
corresponding spaces of $m-$dimensional $l^2-$valued processes
equipped with the norm
\begin{eqnarray}
\|\bm{\phi}_k\|^2_{l^2}&=&\mathbb{E}\left[\int_0^T\sum\limits_{d=1}^{\infty}\sum\limits_{\bm{p}\in\mathbb{N}_d^n}
|\phi_k^{\bm{p}}|^2\right]\qquad k=1,2,\cdots,m,\nonumber\\
\|\bm{\phi}\|_{(l^2)^m}^2&=&\sum\limits_{k=1}^m\|\bm{\phi}_k(t)\|^2_{l^2},\nonumber
\end{eqnarray}
where $\bm{\phi}=(\bm{\phi}_1,\bm{\phi}_2,\cdots,\bm{\phi}_m)'$,
$\bm{\phi}_k=\{\phi_k^{\bm{p}}:\bm{p}\in\mathbb{N}_0^n\}$,
$k=1,2,\cdots,m$ and
$\mathbb{N}_d^n\stackrel{\rm{def}}{=}\{\bm{p}\in\mathbb{N}_0^n:|\bm{p}|=d\}$.
\item [$\bullet$]: Set $\mathcal{H}_T^2=H_T^2\times(M_T^2(l^2))^m$.
\end{itemize}

Following Lin (2011) we introduce power jump monomial processes of
the form
\begin{eqnarray}
X(t)^{(p_1,\cdots,p_n)}\stackrel{\rm{def}}{=}\sum\limits_{0<s\leq
t}(\triangle X_1(s))^{p_1}\cdots(\triangle X_n(s))^{p_n},\nonumber
\end{eqnarray}
The number $|\bm{p}|$ is called the total degree of $X(t)^{\bm{p}}$.
Furthermore define
\begin{eqnarray}
Y(t)^{(p_1,\cdots,p_n)}\stackrel{\rm{def}}{=}X(t)^{(p_1,\cdots,p_n)}-\mathbb{E}[X(t)^{(p_1,\cdots,p_n)}]
=X(t)^{(p_1,\cdots,p_n)}-m_{\bm{p}}t ,\nonumber
\end{eqnarray}
the compensated power jump process of multi-index
$\bm{p}=(p_1,p_2,\cdots,p_n)$ where
$m_{\bm{p}}=\int\prod\limits_{i=1}^n x_i^{p_i}\nu(d\bm{x})$. Under
hypothesis 1, $Y(t)^{(p_1,\cdots,p_n)}$ is a normal martingale,
since for an integrable L\'{e}vy process $Z$, the process
$\{Z_t-E[Z_t], t\geq 0\}$ is a martingale. We call
$Y(t)^{(p_1,\cdots,p_n)}$ the \textsl{Teugels martingale monomial}
of multi-index $(p_1,\cdots,p_n)$.

We can apply the standard Gram-Schmidt process with the graded
lexicographical order to generate a biorthogonal basis
$\{H^{\bm{p}},\bm{p}\in\mathbb{N}^n\}$, such that each
$H^{\bm{p}}(|\bm{p}|=d)$ is a linear combination of the
$Y^{\bm{q}}$, with $|\bm{q}|\leq |\bm{p}|$ and the leading
coefficient equal to $1$. We set
\begin{eqnarray}
H^{\bm{p}}&=&Y^{\bm{p}}+\sum\limits_{\bm{q}\prec\bm{p},|\bm{q}|=|\bm{p}|}c_{\bm{q}}Y^{\bm{q}}
+\sum\limits_{k=1}^{|\bm{p}|-1}\sum\limits_{|\bm{q}|=k}c_{\bm{q}}Y^{\bm{q}},\nonumber
\end{eqnarray}
where $\bm{p}=\{p_1,\cdots,p_n\}$, $\bm{q}=\{q_1,\cdots,q_n\}$ and
$\prec$ represent the relation of graded lexicographical order
between two multi-indexes. Some details about the technique and
theory of orthogonal polynomials of several variables refer to Dunkl
and Xu (2001).

set
\begin{eqnarray}
\textsl{p}(\bm{x})^{\bm{p}}&=&\bm{x}^{\bm{p}}+\sum\limits_{\bm{q}\prec\bm{p},|\bm{q}|=|\bm{p}|}c_{\bm{q}}\bm{x}^{\bm{q}}
+\sum\limits_{k=1}^{|\bm{p}|-1}\sum\limits_{|\bm{q}|=k}c_{\bm{q}}\bm{x}^{\bm{q}},\nonumber\\
\tilde{\textsl{p}}(\bm{x})^{\bm{p}}&=&\bm{x}^{\bm{p}}+\sum\limits_{\bm{q}\prec\bm{p},|\bm{q}|=|\bm{p}|}c_{\bm{q}}\bm{x}^{\bm{q}}
+\sum\limits_{k=2}^{|\bm{p}|-1}\sum\limits_{|\bm{q}|=k}c_{\bm{q}}\bm{x}^{\bm{q}},\nonumber
\end{eqnarray}

Set
\begin{eqnarray}
H^{\bm{p}}(t)&=&\sum\limits_{0<s\leq t}\left((\triangle
X_1)^{p_1}\cdots(\triangle
X_n)^{p_n}+\sum\limits_{\bm{q}\prec\bm{p},|\bm{q}|=|\bm{p}|}c_{\bm{q}}(\triangle
X_1)^{q_1}\cdots(\triangle X_n)^{q_n}\right.\nonumber\\
&&\left.+\sum\limits_{k=1}^{|\bm{p}|-1}\sum\limits_{|\bm{q}|=k}c_{\bm{q}}(\triangle
X_1)^{q_1}\cdots(\triangle
X_n)^{q_n}\right),\nonumber\\
&&-t\mathbb{E}\left[X^{\bm{p}}(1)+\sum\limits_{\bm{q}\prec\bm{p},|\bm{q}|=|\bm{p}|}c_{\bm{q}}X^{\bm{q}}(1)
+\sum\limits_{k=1}^{|\bm{p}|-1}\sum\limits_{|\bm{q}|=k}c_{\bm{q}}X^{\bm{q}}(1)\right]\nonumber\\
&=&\left(c_{\bm{e}_1}X_1(1)+\cdots+c_{\bm{e}_n}X_n(1)\right)+\sum\limits_{0<s\leq
t}\tilde{\textsl{p}}(\triangle
X(s))\nonumber\\
&&-t\mathbb{E}\left[\sum\limits_{0<s\leq
t}\tilde{\textsl{p}}(\triangle
X(s))\right]-t\mathbb{E}\left[c_{\bm{e}_1}X_1(1)+\cdots+c_{\bm{e}_n}X_n(1)\right].\nonumber
\end{eqnarray}

Specially we have
\begin{eqnarray}
H^{\bm{e}_1}(t)&=&c_{\bm{e}_1}(1)(X_1(t)-t\mathbb{E}(X_1(1))),\nonumber\\
H^{\bm{e}_2}(t)&=&c_{\bm{e}_2}(2)(X_2(t)-t\mathbb{E}(X_2(1)))+c_{\bm{e}_1}(2)(X_1(t)-t\mathbb{E}(X_1(1))),\nonumber\\
&\vdots&\\
H^{\bm{e}_n}(t)&=&c_{\bm{e}_n}(n)(X_n(t)-t\mathbb{E}(X_n(1)))+c_{\bm{e}_{n-1}}(n)(X_{n-1}(t)-t\mathbb{E}(X_{n-1}(1)))\nonumber\\
&&+\cdots+c_{\bm{e}_1}(n)(X_1(t)-t\mathbb{E}(X_1(1))).\nonumber
\end{eqnarray}

The main results in Lin (2011) is the Predictable Representation
Property (PRP): Every random variable $F$ in
$L^2(\Omega,\mathscr{F})$ has a representation of the form
\begin{eqnarray}
\begin{array}{rl}
F=&\mathbb{E}(F)+\sum\limits_{d=1}^{\infty}\sum\limits_{\bm{p}\in\mathbb{N}_d^n}
\int_0^T \Phi^{\bm{p}}(s)dH^{\bm{p}}(s)
\end{array}\nonumber
\end{eqnarray}
where $\Phi^{\bm{p}}(s)$ is predictable. It is worthwhile to emphasize that $\Phi^{\bm{p}}(s)$ is not uniqueness, that is to say that $\Phi^{\bm{p}}(s)$ is different for different Gram-Schmidt process, but the form of $\Phi^{\bm{p}}(s)$ is uniqueness. This type of uniqueness is called ``form uniqueness". This result is an extended
version for the corresponding Theorem in Nualart and Schouten
(2000).

\begin{rmk}
If $\nu=0$, we are in the classical
Brownian case and $H^{\bm{p}}(t)=0$, $|\bm{p}|\geq 2$. If $\mu$ has
only mass in $1$, we are in the Poisson case; and also here
$H^{\bm{p}}(t)=0$, $|\bm{p}|\geq 2$. Both case are degenerate cases
in this L\'{e}vy framework.
\end{rmk}

From these observations, it is not so hard to see that the PRP
property shows that financial markets based on a non-Brownian or
non-Poissonian L\'{e}vy process, i.e. with a stock price behaviour
$S_i(t)=S_i(0)exp(rt+X_i(t))$, $i=1,2,\cdots,n$, are so called
incomplete, meaning that perfectly replicating or hedging strategies
do not exists for all relevant contingent claims.

\section{Multidimensional BSDEs Driven by Multidimensional L\'{e}vy Processes}
\label{sec3}
Taking into account the results and notation presented in the
previous section, it seems natural to consider the
\begin{eqnarray}
-d\bm{Y}(t)&=&\bm{f}(t,\bm{Y}(t-),\bm{Z}(t))dt-\sum\limits_{d=1}^{\infty}\sum\limits_{\bm{p}\in\mathbb{N}_d^n}
\bm{z}^{\bm{p}}(s)dH^{\bm{p}}(s), \quad \bm{Y}(T)=\bm{\xi},
\end{eqnarray}
where
\begin{itemize}
\item [$\bullet$]:$\bm{Y}(t)=(Y_1(t),Y_2(t),\cdots,Y_m(t))'$.
\item [$\bullet$]: $\bm{Z}(t)=\{\bm{z}^{\bm{p}}(t)\}_{\bm{p}\in\mathbb{N}_0^n}$, each
component $\bm{z}^{\bm{p}}(t)=(z_1^{\bm{p}},\cdots,z_m^{\bm{p}})'$
is a $m-$variables $\mathscr{F}_t$ predictable function;
\item [$\bullet$] $\bm{f}=(f_1,f_2,\cdots,f_m)':\Omega\times [0,T]\times\mathbb{R}^m\times
  \left(M_T^2(l^2)\right)^m\rightarrow\mathbb{R}^{m}$ is a measurable $m-$dimensional vector function such
  that $\bm{f}(\cdot,\bm{0},\bm{0})\in (H_T^2)^m$.
  \item [$\bullet$] $\bm{f}$ is uniformly Lipschitz in the first two
  components, i.e., there exists $C_k>0$, $k=1,2,\cdots,m$, such that $dt\otimes d\mathbb{P}$
  a.s., for all $(\bm{y}_1,\bm{z}_1)$ and $(\bm{y}_2,\bm{z}_2)$ in $\mathbb{R}^m\times(\bm{l}^2)^m$
  \begin{eqnarray}
   \left|f_k(t,\bm{y}_1,\bm{z}_1)-f_k(t,\bm{y}_2,\bm{z}_2)\right|\leq
   C_k\left(\|\bm{y}_1-\bm{y}_2\|_2+\|\bm{z}_1-\bm{z}_2\|_{(l^2)^m}\right),\qquad k=1,2,\cdots,m.\nonumber
  \end{eqnarray}
  \item [$\bullet$] $\bm{\xi}\in L_T^2(\Omega,\mathbb{P})$.
\end{itemize}

If $(\bm{f},\bm{\xi})$ satisfies the above assumptions, the pair
$(\bm{f},\bm{\xi})$ is said to be \textbf{standard data} for BSDE. A
solution of the BSDE is a pair of processes,
$\{(\bm{Y}(t),\bm{Z}(t)),0\leq t\leq T\}\in H_T^2\times
\left(M_T^2(l^2)\right)^m$ such that the following relation holds
for all $t\in[0,T]$:
\begin{eqnarray}
\bm{Y}(t)=\bm{\xi}+\int_t^T\bm{f}(s,\bm{Y}(s-),\bm{Z}(s))ds-\sum\limits_{d=1}^{\infty}\sum\limits_{\bm{p}\in\mathbb{N}_d^n}
\int_t^T\bm{z}^{\bm{p}}(s)dH^{\bm{p}}(s).
\end{eqnarray}

Note that the progressive measurability of
$\{(\bm{Y}(t),\bm{Z}(t)),0\leq t\leq T\}$ implies that
$(\bm{Y}(0),\bm{Z}(0))$ is deterministic.

When the orthogonal polynomials are fixed, A first key-result
concerns the existence uniqueness of solution of BSDE:
\begin{thm}
Given standard data $(\bm{f},\bm{\xi})$, there exists a unique
form solution $(\bm{Y},\bm{U},\bm{Z})$ which solves the BSDE (3)
\end{thm}
The proof can be found in the Appendix, as the proof of the
continuous dependency of the solution on the final data $\bm{\xi}$
and the function $\bm{f}$.
\begin{thm}
Given standard data $(\bm{f},\bm{\xi})$ and $(\bm{f}',\bm{\xi}')$,
let $(\bm{Y},\bm{Z})$ and $(\bm{Y}',\bm{Z}')$ be the unique adapted
form solutions of the BSDE(3) corresponding to $(\bm{f},\bm{\xi})$ and
$(\bm{f}',\bm{\xi}')$. Then
\begin{eqnarray}
\mathbb{E}\left[\int_0^T\left(\|\bm{Y}(s-)-\bar{\bm{Y}}(s-)\|^2+\sum\limits_{d=1}^{\infty}\sum\limits_{\bm{p}\in\mathbb{N}_d^n}
\|\bm{z}^{\bm{p}}(s)-\bar{\bm{z}}^{\bm{p}}(s)\|^2\right)ds\right]\nonumber\\
\leq C\left(\mathbb{E}[\|\bm{\xi}-\bar{\bm{\xi}}\|^2]+
\mathbb{E}\left[\int_0^T\|\bm{f}(s,\bm{Y}(s-),\bm{Z}(s))-\bar{\bm{f}}(s,\bm{Y}(s-),\bm{Z}(s))\|^2ds\right]\right).\nonumber
\end{eqnarray}
\end{thm}

The definition of ``form solution" in Theorem 1 and Theorem 2 is deferred until later in the proofs of Theorem 1 and Theorem 2 in Appendix.

\section{Applications}
\label{sec4}

Suppose our n-dimensional L\'{e}vy process $\bm{X}(t)$ has no
Brownian part, i.e. $\bm{X}(t)=\bm{a}t+\bm{L}(t)$, where
$\bm{a}=(a_1,\cdots,a_n)'$ and $\bm{L}(t)$ is n-dimensional pure
jump process with L\'{e}vy measure $v(d\bm{x})$.
\subsection{Clark-Ocone Formula and Feynman-Kac Formula}
Let us consider the simple case of a BSDE where $\bm{f}=\bm{0}$, and
the terminal random vector $\bm{\xi}$ is a function of $\bm{X}(T)$,
that is,
\begin{eqnarray}
-d\bm{Y}(t)=-\sum\limits_{d=1}^{\infty}\sum\limits_{\bm{p}\in\mathbb{N}_d^n}\bm{z}^{\bm{p}}(t)dH^{\bm{p}}(t),\qquad
\bm{Y}(T)=\bm{g}(\bm{X}(T))\nonumber
\end{eqnarray}
or equivalently
\begin{eqnarray}
\bm{Y}(t)=\bm{g}(\bm{X}(T))-\sum\limits_{d=1}^{\infty}\sum\limits_{\bm{p}\in\mathbb{N}_d^n}
\int_t^T\bm{z}^{\bm{p}}(s)dH^{\bm{p}}(s),
\end{eqnarray}
where $\bm{g}=(g_1,g_2,\cdots,g_m)'$ and
$\mathbb{E}(\|\bm{g}(\bm{X}(T))\|^2)<\infty$. Let
$\theta_k=\theta_k(t,\bm{x})$, $k=1,2,\cdots,m$, be the solution of
the following PDIE(Partial Differential Integral Equation) with
terminal value $\bm{g}_k$:
\begin{eqnarray}
\frac{\partial\theta_k}{\partial
t}(t,\bm{x})+\int_{\mathbb{R}^n}\left(\theta_k(t,\bm{x}+\bm{y})-\theta_k(t,\bm{x})
-\frac{\partial\theta_k}{\partial\bm{x}}(t,\bm{x})\cdot\bm{y}\right)\nu(d\bm{y})+\tilde{\bm{a}}\cdot\frac{\partial\theta_k}
{\partial\bm{x}}(t,\bm{x})=0,\nonumber\\
\theta_k(T,\bm{x})=g_k(\bm{x}),
\end{eqnarray}
where $\tilde{\bm{a}}=\bm{a}+\int_{\{\|\bm{y}\|\geq
1\}}\bm{y}\nu(d\bm{y})$. Set
\begin{eqnarray}
\theta_k^{(1)}(t,\bm{x},\bm{y})=\theta_k(t,\bm{x}+\bm{y})-\theta_k(t,\bm{x})-\frac{\partial\theta_k}{\partial\bm{x}}(t,\bm{x})\cdot\bm{y}.
\end{eqnarray}
The following result is a version of the Clark-Ocone formula for
functions of a L\'{e}vy process. Again the proof can be found in the
Appendix.
\begin{pro}
Suppose that $\theta_k$ is a $C^{1,2}$ function for $k=1,2,\cdots,m$
such that $\frac{\partial\theta_k}{\partial\bm{x}}$ and
$\frac{\partial^2\theta_k}{\partial\bm{x}^2}$ are bounded by a
polynomial function of $\bm{x}$, uniformly in $t$, then the unique
adapted form solution of (4) is given by
\begin{eqnarray}
\begin{array}{lcl}
  Y_k(t) & = & \theta_k(t,\bm{X}(t)) \\
  z_k^{\bm{p}} & = & \int_{\mathbb{R}^n}\theta_k^{(1)}(t,\bm{X}(t-),\bm{y})
  \textsl{p}^{\bm{p}}(\bm{y})\nu(d\bm{y}),\qquad |\bm{p}|\geq 2,\\
  z_k^{\bm{e}_i} &=&\int_{\mathbb{R}^n}\theta_k^{(1)}(t,\bm{X}(t-),\bm{y})
  \textsl{p}^{\bm{e}_i}(\bm{y})\nu(d\bm{y})+\sum\limits_{j=1}^n\frac{\partial\theta_k}{\partial
  x_j}\tilde{c}_{ij},\quad i=1,2,\cdots,n
\end{array}\nonumber
\end{eqnarray}
where $\theta_k=\theta_k(t,\bm{x})$ for $k=1,2,\cdots,m$ are the
solutions of the system of PDIE(5),
$\theta_k^{(1)}(t,\bm{x},\bm{y})$ for $k=1,2,\cdots, m$ are given by
(6) and $[\tilde{c}_{ij}]$ is the inverse matrix of the coefficient
matrix Gram-Schmidt transformation (1).
\end{pro}
Now by taking expectations we derive that the solution
$\bm{\theta}(t,\bm{x})=(\theta_1(t,\bm{x}),\cdots,\theta_m(t,\bm{x}))'$
to our PDIE(5) equation has the stochastic representation
\begin{eqnarray}
\bm{\theta}(t,\bm{x})=\mathbb{E}\left[\bm{g}(\bm{X}(T))|\bm{X}(t)=\bm{x}\right].\nonumber
\end{eqnarray}
This is an extension of the classical Feynman-Kac Formula.

If $\int_{\mathbb{R}^n}\|\bm{y}\|\nu(d\bm{y})<\infty$, and we take
$\bm{a}=\int_{\{\|\bm{y}\|<1\}}\bm{y}\nu(d\bm{y})$, then the
equation (5) reduces to
\begin{eqnarray}
\frac{\partial\bm{\theta}}{\partial
t}(t,\bm{x})+\int_{\mathbb{R}^n}\left(\bm{\theta}(t,\bm{x}+\bm{y})-\bm{\theta}(t,\bm{x})\right)\nu(d\bm{y})=0,\nonumber\\
\bm{\theta}(T,\bm{x})=\bm{g}(\bm{x}),\nonumber
\end{eqnarray}
and we have for $k=1,2,\cdots,m$
\begin{eqnarray}
\begin{array}{lcl}
  z_k^{\bm{p}} & = &
  \int_{\mathbb{R}^n}(\theta_k(t,\bm{X}(t-)+\bm{y})-\theta_k(t,\bm{X}(t-)))
  \textsl{p}^{\bm{p}}(\bm{y})\nu(d\bm{y}),\qquad |\bm{p}|\geq 2,\\
  z_k^{\bm{e}_i}
  &=&\int_{\mathbb{R}^n}(\theta_k(t,\bm{X}(t-)+\bm{y})-\theta_k(t,\bm{X}(t-)))
  \textsl{p}^{\bm{e}_i}(\bm{y})\nu(d\bm{y})+\sum\limits_{j=1}^n\frac{\partial\theta_k}{\partial
  x_j}\tilde{c}_{ij},\quad i=1,2,\cdots,n.
\end{array}\nonumber
\end{eqnarray}

\noindent\textbf{Example:} In this example, we define a
two-dimensional Poisson process by using L\'{e}vy copulas. All the
concepts and notations are adopted from the Kallsen and Tankov
(2006). Here the two marginal processes are given respectively by
two Poisson process $N_i(t)$ where $i=1,2$. In particular, the
L\'{e}vy copula
$F(u_1,u_2):\bar{\mathbb{R}}^n\rightarrow\bar{\mathbb{R}}$ is taken
as
\begin{eqnarray}
F(u_1,u_2)=\left(|u_1|^{-\mu}+|u_2|^{-\mu}\right)^{-1/\mu}(\eta
I_{\{u_1u_2\geq 0\}}-(1-\eta)I_{\{u_1u_2<0\}})
\end{eqnarray}
where $\mu>0$ and any $\eta\in [0,1]$. It defines a two parameter
family of L\'{e}vy copulas which resembles the Clayton family of
ordinary copulas. Thus its L\'{e}vy measure can be calculated as:
\begin{eqnarray}
\nu(dx_1dx_2)=\frac{\eta(1+\mu)(\lambda_1\lambda_2)^\mu}
{\left(\lambda_1^\mu+\lambda_2^\mu\right)^{\frac{1}{\mu}+2}}\delta_1(x_1)\delta_1(x_2).
\end{eqnarray}

Set two compensated Poisson process $X_i(t)=N_i(t)-\lambda_it$ for
$i=1,2$, then we can construct a set of martingales $H^{\bm{p}}$,
$\bm{p}\in\mathbb{N}_0^n$, $i=1,2$ by orthogonalizing procedure
proposed by Lin (2011) such that $H^{\bm{p}}$ is strongly orthogonal
to $H^{\bm{q}}$ with respect to L\'{e}vy measure $\nu(dx_1dx_2)$,
for $\bm{p}\neq\bm{q}$. Moreover the PDIE (5) reduces to
\begin{eqnarray}
\frac{\partial\bm{\theta}}{\partial
t}(t,\bm{x})+\left(\bm{\theta}(t,\bm{x}+\bm{1})-\bm{\theta}(t,\bm{x})
\right)-\frac{\partial\bm{\theta}}{\partial\bm{x}}(t,\bm{x})\cdot\bm{1}\frac{\eta(1+\mu)(\lambda_1\lambda_2)^\mu}
{\left(\lambda_1^\mu+\lambda_2^\mu\right)^{\frac{1}{\mu}+2}}=0,\nonumber\\
\bm{\theta}(T,\bm{x})=\bm{g}(\bm{x}),\nonumber
\end{eqnarray}
The Clark-Ocone Formula is now given by
\begin{eqnarray}
g_k(\bm{X}(T))&=&\mathbb{E}[g_k(\bm{X}(T))]+\sum_{i=1}^n
\int_t^T\left(\theta_k(s,\bm{X}(s)+\bm{1})-\theta_k(s,\bm{X}(s))
\right)d\bm{X}_i(s) \nonumber
\end{eqnarray}

\subsection{Nonlinear Clark-Haussman-Ocone Formula and Feynman-Kac Formula}
Let us consider the BSDE
\begin{eqnarray}
-d\bm{Y}(t)&=&\bm{f}(t,\bm{Y}(t),\bm{Z}(t))-\sum\limits_{d=1}^{\infty}\sum\limits_{\bm{p}\in\mathbb{N}_d^n}
\bm{z}^{\bm{p}}(t)dH^{\bm{p}}(t), \quad \bm{Y}(T)=\bm{g}(\bm{X}(T))
\end{eqnarray}
or equivalently
\begin{eqnarray}
\bm{Y}(t)=\bm{g}(\bm{X}(T))+\int_t^T\bm{f}(s,\bm{Y}(s-),\bm{Z}(s))ds-
\sum\limits_{d=1}^{\infty}\sum\limits_{\bm{p}\in\mathbb{N}_d^n}
\int_t^T\bm{z}^{\bm{p}}(s)dH^{\bm{p}}(s),\nonumber
\end{eqnarray}
Suppose that $\theta_k=\theta_k(t,\bm{x})$ for $k=1,2,\cdots,m$
satisfy the following system of PDIE:
\begin{eqnarray}
\frac{\partial\theta_k}{\partial
t}(t,\bm{x})+\int_{\mathbb{R}^n}\theta_k^{(1)}(t,\bm{x},\bm{y})\nu(d\bm{y})+\tilde{\bm{a}}\cdot\frac{\partial\theta_k}
{\partial\bm{x}}(t,\bm{x})+f_k\left(t,\theta_k(t,\bm{x}),\bm{\Theta}_k(t,\bm{x})\right)=0,\\
\theta_k(T,\bm{x})=g_k(\bm{x}).\nonumber
\end{eqnarray}
where as in the previous section, we define
$\theta_k^{(1)}(t,\bm{x},\bm{y})$ by (6), and
$\Theta_k(t)=\{\theta_k^{\bm{p}}(t)\}_{\bm{p}\in\mathbb{N}_0^n}$,
where
\begin{eqnarray}
\theta_k^{\bm{p}}(t,\bm{x})&=&\int_{\mathbb{R}^n}\theta_k^{(1)}(t,\bm{x},\bm{y})\textsl{p}^{\bm{p}}(\bm{y})\nu(d\bm{y}),
\end{eqnarray}

\begin{pro}
Suppose that $\theta_k$ is a $C^{1,2}$ function such that
$\frac{\partial\theta_k}{\partial\bm{x}}$ and
$\frac{\partial^2\theta_k}{\partial\bm{x}^2}$ are bounded by a
polynomial function of $\bm{x}$, uniformly in $t$, then the unique
adapted form solution of (9) is given by
\begin{eqnarray}
\begin{array}{lcl}
  Y_k(t) & = & \theta_k(t,\bm{X}(t)) \\
  z_k^{\bm{p}} & = & \int_{\mathbb{R}^n}\theta_k^{(1)}(t,\bm{X}(t-),\bm{y})
  \textsl{p}^{\bm{p}}(\bm{y})\nu(d\bm{y}),\qquad |\bm{p}|\geq 2, \\
  z_k^{\bm{e}_i} &=&\int_{\mathbb{R}^n}\theta_k^{(1)}(t,\bm{X}(t-),\bm{y})
  \textsl{p}^{\bm{e}_i}(\bm{y})\nu(d\bm{y})+\sum\limits_{j=1}^n\frac{\partial\theta}{\partial
  x_j}\tilde{c}_{ij},\quad i=1,2,\cdots,n.
\end{array}\nonumber
\end{eqnarray}
where $\theta_k=\theta_k(t,\bm{x})$ for $k=1,2,\cdots,m$ are the
solution of the system of PDIE(10),
$\theta_k^{(1)}(t,\bm{x},\bm{y})$ is given by (6) and
$[\tilde{c}_{ij}]$ is the inverse matrix of the coefficient matrix
Gram-Schmidt transformation (1).
\end{pro}
Notice that taking expectations we get

\begin{eqnarray}
\theta_k(t,\bm{x})=\mathbb{E}\left[g_k(\bm{X}(T))|\bm{X}(t)=\bm{x}\right]
+\mathbb{E}\left[\int_t^Tf_k\left(s,\theta_k(s,\bm{X}(s-)),
\bm{\Theta}_k(s,\bm{X}(s-))\right)ds|\bm{X}(t)=\bm{x}\right].\nonumber
\end{eqnarray}

\noindent\textbf{Example:} Consider again the very special case
where have a two-dimensional Poisson process $(N_1(t),N_2(t))$ by
using L\'{e}vy copulas (7). Set $X_i(t)=N_i(t)-\lambda_it$ where
$i=1,2$. Then the PDIE (10) reduces to
\begin{eqnarray}
\frac{\partial\theta_k}{\partial
t}(t,\bm{x})+\left(\theta_k(t,\bm{x}+\bm{1})-\theta_k(t,\bm{x})
-\frac{\partial\theta_k}{\partial\bm{x}}(t,\bm{x})\cdot\bm{1}\right)\frac{\eta(1+\mu)(\lambda_1\lambda_2)^\mu}
{\left(\lambda_1^\mu+\lambda_2^\mu\right)^{\frac{1}{\mu}+2}}\nonumber\\
+f_k\left(t,\theta_k(t,\bm{x}),\theta_k(t,\bm{x}+\bm{1})-\theta_k(t,\bm{x})\right)=0\nonumber\\
\theta_k(T,\bm{x})=g_k(\bm{x}),\nonumber\\
k=1,2,\cdots,m\nonumber
\end{eqnarray}
And we derive the nonlinear Feynman-Kac Formula:
\begin{eqnarray}
\theta_k(t,\bm{x})=\mathbb{E}\left[g_k(\bm{X}(T))|\bm{X}(t)=\bm{x}\right]
+\mathbb{E}\left[\int_t^Tf_k\left(s,\theta_k(s,\bm{X}(s-)),
\bm{\Theta}_k(s,\bm{X}(s-))\right)ds|\bm{X}(t)=\bm{x}\right].\nonumber
\end{eqnarray}
where $k=1,2,\cdots,m$.

\subsection{Option Pricing}
In the last two decades several particular choices for
one-dimensional non-Brownian L\'{e}vy processes were proposed. Madan
and Seneta (1990) have proposed a L\'{e}vy process with variance
gamma distributed increments. We mention also the Hyperbolic Model
proposed by Eberlein and Keller (1995). In the same year
Barndorff-Nielsen (1995) proposed the normal inverse Gaussian
L\'{e}vy process. The CMGY model was also introduced in Carr et al.
(2000). Finally, we mention the Meixner model (see Grigelionis
(1999) and Schoutens (2001)). All models give a much better fit to
the data and lead to an improvement with respect to the
Black-Scholes model.

Multidimensional models with jumps are more difficult to construct
than one-dimensional ones. A simple method to introduce jumps into a
multidimensional model is to take a multivariate Brownian motion and
time change it with a univariate subordinator (refer to Cont and
Tankov (2004)). The multidimensional versions of the models include
variance gamma, normal inverse Gaussian and generalized hyperbolic
processes. The principal advantage of this method is its simplicity
and analytic tractability; in particular, processes of this type are
easy to simulate. Another method to introduce jumps into a
multidimensional model is so-called method of L\'{e}vy copulas
proposed by Kallsen and Tankov (2006). The principle advantage in
this way lies in that the dependence among components of the
multidimensional L\'{e}vy processes can be completely characterized
with a L\'{e}vy copula. This allows us to give a systematic method
to construct multidimensional L\'{e}vy processes with specified
dependence.

Here we define a multivariate Meixner process by using L\'{e}vy
copulas, and all the concepts and notations are adopted from the
Kallsen and Tankov (2006). In particular, for $n\geq 2$, the
L\'{e}vy copula $F(u_1,\cdots,u_n):
\bar{\mathbb{R}}^n\rightarrow\bar{\mathbb{R}}$ is taken as
\begin{eqnarray}
F(u_1,\cdots,u_n)=2^{2-n}\left(\sum\limits_{j=1}^n|u_j|^{-\mu}\right)^{-1/\mu}(\eta
I_{\{u_1\cdots u_n\geq 0\}}-(1-\eta)I_{\{u_1\cdots u_n< 0\}})
.\nonumber
\end{eqnarray}
It defines a two parameter family of L\'{e}vy copulas which
resembles the Clayton family of ordinary copulas. It is in fact a
L\'{e}vy copula homogeneous of order 1, for any $\mu>0$ and any
$\eta\in [0,1]$.

In addition, we know that if the tail integrals $U_i(x_i)$ related
to the marginal L\'{e}vy densities $\nu_i(dx_i)$, $i=1,\cdots,n$,
are absolutely continuous, we can compute the L\'{e}vy density of
the L\'{e}vy copula process by differentiation as follows:
\begin{eqnarray}
\nu(dx_1,\cdots,dx_n)=\partial_1\cdots\partial_nF|_{\xi_1=U_1(x_1),\cdots,\xi_n=U_n(x_n)}\nu_1(dx_1)\cdots\nu_1(dx_n)\nonumber
\end{eqnarray}
where $\nu_1(dx_1),\cdots,\nu_n(x_n)$ are marginal L\'{e}vy
densities.

Hence we are able to construct a $n-$variate L\'{e}vy process
$\bm{X}(t)=(X_1(t),X_2(t),\cdots,X_n(t))'$ with state space
$\mathbb{R}^n$ and characteristic triple
$(\mathbb{P},\nu_{\mathbb{P}},\bm{a})$ by using the above
\textsl{complete dependent L\'{e}vy copula or Clayton L\'{e}vy
copulas} proposed by Kallsen and Tankov (2006). Here, in particular,
we will consider L\'{e}vy marginal measure $\nu_i$ of Meixner type
for L\'{e}vy copula. We recall the marginal density $f_i$, the
cumulative generating function $K_i$, the drift $a_i$, and the
marginal L\'{e}vy measure $\nu_i$, for the ith component of the
Meixner Process $\{X_i(t), t\geq 0\}$, $i=1,2,\cdots,n$, are
illustrated as follows,
\begin{eqnarray}
\frac{\mathbb{P}_{Meix}(dx_i)}{dx_i}&=&f_i(x_i;\alpha_i,\beta_i,\delta_i,\mu_i)
=\frac{\left(2cos\frac{\beta_i}{2}\right)^{2\delta_i}e^{\frac{\beta_i(x_i-\mu_i)}{\alpha_i}}|\Gamma\left(\delta_i+\frac{i(x_i-\mu_i)}{\alpha_i}\right)|^2}
{\pi\alpha_i\Gamma(2\delta_i)},\nonumber\\
K_i(\theta_i;\alpha_i,\beta_i,\delta_i,\mu_i)&=&\mu_i\theta_i+2\delta_i\left(log
cos\frac{\beta_i}{2}-log
cos\frac{\alpha_i\theta_i+\beta_i}{2}\right),\nonumber\\
a_i(\alpha_i,\beta_i,\delta_i,\mu_i)&=&\mu_i+\alpha_i\delta_itan\frac{\beta_i}{2}-2\delta_i\int_1^{\infty}
\frac{sinh\frac{\beta_ix_i}{\alpha_i}}{sinh\frac{\pi
x_i}{\alpha_i}}dx_i,\nonumber\\
\nu_i(dx_i;\alpha_i,\beta_i,\delta_i,\mu_i)&=&\frac{\delta_ie^{\frac{\beta_ix_i}{\alpha_i}}}{x_isinh\frac{\pi
x_i}{\alpha_i}}dx_i,\nonumber
\end{eqnarray}
where $\alpha_i>0$, $-\pi<\beta_i<\pi$, $\mu_i\in\mathbb{R}$, and
$\delta_i>0$.

From the form of the cumulative generating function one easily
deduces that the density at any time $t$ can be calculated by
multiplying the parameters $\delta_i$ and $\mu_i$ by $t$ for both
cases.

In addition, Reich, etc.(2009) had proved that $e^{X_i}$ is a
martingale with respect to the canonical filtration $\mathscr{F}$ of
$\bm{X}$ if only if
\begin{eqnarray}
\frac{\sigma_{ij}}{2}+a_j+\int_{\mathbb{R}^n}(e^{z_j}-1-z_jI_{|\bm{z}|\leq
1})\nu(d\bm{z})=0.\nonumber
\end{eqnarray}

Following Reich, etc.(2009), we assume a risk-neutral market which
consists of one riskless asset (the bond) with price process given
by $B(t)=e^{rt}$, where $r$ is compound interest rate, and $n(\geq
1)$ risky assets (the stocks), with price process:
\begin{eqnarray}
S_i(t)=S_i(0)exp(rt+X_i(t)),\quad i=1,2,\cdots,n, \nonumber
\end{eqnarray}
where $\bm{X}(t)=(X_1(t),X_2(t),\cdots,X_n(t))'$ is a $n-$variate
L\'{e}vy process and characteristic triple
$(\mathbb{P},\nu_{\mathbb{P}},\bm{a})$ under a risk-neutral measure
$\mathbb{P}$ such that $e^{X_i}$ is a martingale with respect to the
canonical filtration
$\mathscr{F}_t^0\stackrel{\rm{def}}{=}\sigma(\bm{X}(s),s\leq t)$,
$t\geq 0$, of the multivariate process $\bm{X}$. Denote by
$\mathbb{P}(d\bm{x})$ the probability measure of $\bm{X}(1)$.

We consider a European option with maturity $T<\infty$ and payoff
$G(\bm{S})$ where $\bm{S}=(S_1,\cdots,S_n)'$, and we assume that
$G(\bm{S})$ is Lipschitz. According to the fundamental theorem of
asset pricing (see Delbaen and Schachermayer (1994)) the value
$V(t,\bm{S})$ of this option is given by
$$V(t,\bm{S})=\mathbb{E}[e^{-r(T-t)}G(\bm{S}(T))|\bm{S}(t)=\bm{S}].$$
If $V(t,\bm{S})$ satisfies
$$V(t,\bm{S})\in C^{1,2}\left((0,T)\times\mathbb{R}_{>0}^n\right)\cap C^0\left([0,T]\times\mathbb{R}_{\geq 0}^n\right)$$
Then Reich,etc.(2009) had proved that $V(t,\bm{S})$ is a classical
solution of the backward Kolmogorov equation:
\begin{eqnarray}
\frac{\partial V}{\partial
t}(t,\bm{S})+\frac{1}{2}\sum\limits_{i,j=1}^nS_iS_j\sigma_{ij}\frac{\partial^2V}{\partial
S_i\partial S_j}(t,\bm{S})+r\sum\limits_{i=1}^nS_i\frac{\partial
V}{\partial
S_i}(t,\bm{S})-r V(t,\bm{S})\nonumber\\
+\int_{\mathbb{R}^n}\left(V(t,\bm{S}\cdot
e^{\bm{z}})-V(t,\bm{S})-\sum\limits_{i=1}^nS_i(e^{z_i}-1)\frac{\partial
V}{\partial S_i}(t,S)\right)\nu(d\bm{z})=0,\nonumber
\end{eqnarray}
on $(0,T)\times\mathbb{R}_{\geq 0}^n$ where $V(t,\bm{S}\cdot
e^{\bm{z}})\stackrel{\rm{def}}{=}V(t,S_1e^{z_1},\cdots,S_ne^{z_n})$,
and the terminal condition is given by
\begin{eqnarray}
V(T,\bm{S})=g(\bm{S}),\qquad \forall \bm{S}\in\mathbb{R}_{\geq
0}^n.\nonumber
\end{eqnarray}

\section{Appendix: Proofs of the Results}
\noindent\textbf{Proof of Theorem 1:}

We define a mapping $\bm{\Phi}$ from $\mathcal{H}_T^2$ into itself
such that $(\bm{Y},\bm{Z})\in\mathcal{H}_T^2$ is a solution of the
BSDE if only if it is a fixed point of $\bm{\Phi}$. Given
$(\bm{U},\bm{V})\in\mathcal{H}_T^2$, we define
$(\bm{Y},\bm{Z})=\bm{\Phi}(\bm{U},\bm{V})$ as follows:
\begin{eqnarray}
\bm{Y}(t)=E\left[\bm{\xi}+\int_t^T\bm{f}(s,\bm{U}(s-),\bm{V}(s))ds|\mathcal{F}_t\right],\quad
0\leq t\leq T,\nonumber
\end{eqnarray}
and $\{\bm{Z}(t), 0\leq t\leq T\}$ is given by the martingale
representation of Lin (2011) applied to the square integrable random
variable
\begin{eqnarray}
\bm{\xi}+\int_0^T\bm{f}(s,\bm{U}(s-),\bm{V}(s))ds,\nonumber
\end{eqnarray}
i.e.,
\begin{eqnarray}
\begin{array}{rl}
&\bm{\xi}+\int_0^T\bm{f}(s,\bm{U}(s-),\bm{V}(s))ds\nonumber\\
=&E\left[\bm{\xi}+\int_0^T\bm{f}(s,\bm{U}(s-),\bm{V}(s))ds\right]
+\sum\limits_{d=1}^{\infty}\sum\limits_{\bm{p}\in\mathbb{N}_d^n}
\int_0^T \phi^{\bm{p}}(s)dH^{\bm{p}}(s),\nonumber
\end{array}
\end{eqnarray}

In according to the result in Lin (2011), a different Gram-Schmidt process will generate a different $\Phi^{\bm{p}}(s)$. Here and hereafter, for a given $(\bm{U},\bm{V})\in\mathcal{H}_T^2$, we always mean that the $\Phi^{\bm{p}}(s)$ is generated by arbitrarily selecting a fixed Gram-Schmidt process. It is to say that the Gram-Schmidt process is fixed once the Gram-Schmidt process is first arbitrarily selected. In the following, the existence and uniqueness of the solution for BSDE(3) will related to this selected $\Phi^{\bm{p}}(s)$. Although a $\Phi^{\bm{p}}(s)$ is different for a different Gram-Schmidt process, the expression forms of $\Phi^{\bm{p}}(s)$ are the same, the corresponding solutions are called ``form solution".

Taking the conditional expectation with respect to $\mathcal{F}_t$
in the last identity yields
\begin{eqnarray}
\bm{Y}_t+\int_0^t\bm{f}(s,\bm{U}(s-),\bm{V}(s))ds=\bm{Y}_0+\sum\limits_{d=1}^{\infty}\sum\limits_{\bm{p}\in\mathbb{N}_d^n}
\int_0^t\phi^{\bm{p}}(s)dH^{\bm{p}}(s),\nonumber
\end{eqnarray}
from which we deduce that
\begin{eqnarray}
\bm{Y}_t=\bm{\xi}+\int_t^T\bm{f}(s,\bm{U}(s-),\bm{V}(s))ds-\sum\limits_{d=1}^{\infty}\sum\limits_{\bm{p}\in\mathbb{N}_d^n}
\int_t^T \phi^{\bm{p}}(s)dH^{\bm{p}}(s),\nonumber
\end{eqnarray}
and we have shown that $(\bm{Y},\bm{Z})\in\mathcal{H}_T^2$ solves
our BSDE if only if it is a fixed point of $\Phi$.

Next we prove that $\Phi$ is a strict contraction on
$\mathcal{H}_T^2$ equipped with the norm
\begin{eqnarray}
\|(\bm{Y},\bm{Z})\|_\beta=\left(\int_0^Te^{\beta
s}\left(\|\bm{Y}(s-)\|^2+\|\bm{Z}(s)\|^2\right)ds\right)^{1/2},\nonumber
\end{eqnarray}
for a suitable $\beta>0$. Let $(\bm{U},\bm{V})$ and
$(\bm{U}',\bm{V}')$ be two elements of $\mathcal{H}_T^2$ and set
$\Phi(\bm{U},\bm{V})=(\bm{Y},\bm{Z})$ and
$\Phi(\bm{U}',\bm{V}')=(\bm{Y}',\bm{Z}')$. Denote
$(\overline{\bm{U}},\overline{\bm{V}})=(\bm{U}-\bm{U}',\bm{V}-\bm{V}')$
and
$(\overline{\bm{Y}},\overline{\bm{Z}})=(\bm{Y}-\bm{Y}',\bm{Z}-\bm{Z}')$.

Applying It$\hat{o}$'s formula from $s=t$ to $s=T$, to $e^{\beta
s}\|\bm{Y}(s)-\bm{Y}(s)'\|^2$, it follows that
\begin{eqnarray}
e^{\beta t}\|\bm{Y}(t)-\bm{Y}(t)'\|^2&=&-\beta\int_t^Te^{\beta
s}\|\bm{Y}(s-)-\bm{Y}'(s-)\|^2ds\nonumber\\
&&-2\int_t^Te^{\beta s}(\bm{Y}(s-)-\bm{Y}'(s-))\cdot d(\bm{Y}(s)-\bm{Y}(s)')\nonumber\\
&&-\sum\limits_{j=1}^m\int_t^Te^{\beta
s}d[\bm{Y}_j-\bm{Y}_j',\bm{Y}_j-\bm{Y}_j'](s).
\end{eqnarray}
We have
\begin{eqnarray}
-d(\bm{Y}(t)-\bm{Y}(t)')&=&(\bm{f}(t,\bm{U}(t-),\bm{V}(t))-\bm{f}(t,\bm{U}(t-)',\bm{V}(t)'))dt\nonumber\\
&&+\sum\limits_{d=1}^{\infty}\sum\limits_{\bm{p}\in\mathbb{N}_d^n}
\overline{\bm{\phi}}^{\bm{p}}(t)dH^{\bm{p}}(t),\nonumber
\end{eqnarray}
\begin{eqnarray}
\begin{array}{rl}
&d[\bm{Y}_j-\bm{Y}_j',\bm{Y}_j-\bm{Y}_j'](t)\nonumber\\
=&\sum\limits_{d=1}^{\infty}\sum\limits_{\bm{p}\in\mathbb{N}_d^n}\sum\limits_{e=1}^{\infty}\sum\limits_{\bm{q}\in\mathbb{N}_e^n}
\overline{\bm{\phi}}_j^{\bm{p}}(t)\overline{\bm{\phi}}_j^{\bm{q}}(t)d[H^{\bm{p}},H^{\bm{q}}](t),\qquad
j=1,2,\cdots,m\nonumber
\end{array}
\end{eqnarray}
where the symbol ``$\circ$" represents the Hadamard-Schur product
for two vectors.
\begin{eqnarray}
<H^{\bm{p}},H^{\bm{q}}>(t)=\delta_{\bm{p},\bm{q}}t. \nonumber
\end{eqnarray}
Hence, taking expectations in (12), we have
\begin{eqnarray}
\begin{array}{rl}
&\mathbb{E}\left[e^{\beta t}\|\bm{Y}(t)-\bm{Y}(t)'\|^2\right]
+\sum\limits_{j=1}^m\sum\limits_{d=1}^{\infty}\sum\limits_{\bm{p}\in\mathbb{N}_d^n}\mathbb{E}\left[\int_t^Te^{\beta
s}\overline{\bm{\phi}}_j^{\bm{p}}(s)^2ds\right]\nonumber\\
=&-\beta\mathbb{E}\left[\int_t^Te^{\beta
s}\|\bm{Y}(s-)-\bm{Y}'(s-)\|^2ds\right]\nonumber\\
&+2\mathbb{E}\left[\int_t^Te^{\beta
s}(\bm{Y}(s-)-\bm{Y}'(s-))\cdot(\bm{f}(s,\bm{U}(s-),\bm{V}(s))-\bm{f}(s,\bm{U}(s-)',\bm{V}(s)'))ds\right].\nonumber
\end{array}
\end{eqnarray}
Using the fact that $\bm{f}$ is Lipschitz with constant $C$ yields
\begin{eqnarray}
\begin{array}{rl}
&\mathbb{E}\left[e^{\beta t}\|\bm{Y}(t)-\bm{Y}(t)'\|^2\right] +
\mathbb{E}\left[\int_t^Te^{\beta s}\|\overline{\bm{\phi}}(s)\|^2ds\right]\nonumber\\
\leq&-\beta\mathbb{E}\left[\int_t^Te^{\beta
s}\|\bm{Y}(s-)-\bm{Y}'(s-)\|^2ds\right]\nonumber\\
&+2C\mathbb{E}\left[\int_t^Te^{\beta
s}\|\bm{Y}(s-)-\bm{Y}'(s-)\|\cdot\left(|\bm{U}(s-)-\bm{U}'(s-)|+\|\overline{\bm{V}}(s)\|\right)ds\right].\nonumber
\end{array}
\end{eqnarray}
If we now use the fact that for every $c>0$ and $a,b\in\mathbb{R}$
we have that $2ab\leq ca^2+\frac{1}{c}b^2$ and $(a+b)^2\leq
2a^2+2b^2$, we obtain
\begin{eqnarray}
\begin{array}{rl}
&\mathbb{E}\left[e^{\beta t}\|\bm{Y}(t)-\bm{Y}(t)'\|^2\right]
+\mathbb{E}\left[\int_t^Te^{\beta s}\|\overline{\bm{\phi}}(s)\|^2ds\right]\nonumber\\
\leq&(4C^2-\beta)\mathbb{E}\left[\int_t^Te^{\beta
s}\|\bm{Y}(s-)-\bm{Y}'(s-)\|^2ds\right]\nonumber\\
&+\frac{1}{2}\mathbb{E}\left[\int_t^Te^{\beta
s}\left(|\bm{U}(s-)-\bm{U}'(s-)|^2+\|\overline{\bm{V}}(s)\|^2\right)ds\right].\nonumber
\end{array}
\end{eqnarray}
Taking now $\beta=4C^2+1$, and noting that $e^{\beta
t}\mathbb{E}\left[\|\bm{Y}(t)-\bm{Y}(t)'\|^2\right]\geq 0$, we
finally derive
\begin{eqnarray}
\begin{array}{rl}
&\mathbb{E}\left[\int_t^Te^{\beta
s}\|\bm{Y}(s)-\bm{Y}(s)'\|^2ds\right]
+\mathbb{E}\left[\int_t^Te^{\beta s}\|\overline{\bm{\phi}}(s)\|^2ds\right]\nonumber\\
\leq&\frac{1}{2}\mathbb{E}\left[\int_t^Te^{\beta
s}\left(|\bm{U}(s-)-\bm{U}'(s-)|^2+\|\overline{\bm{V}}(s)\|^2\right)ds\right].\nonumber
\end{array}
\end{eqnarray}
that is
\begin{eqnarray}
\|(\bm{Y},\bm{Z})\|_\beta^2\leq\frac{1}{2}\|(\bm{U},\bm{V})\|_\beta^2,\nonumber
\end{eqnarray}
from which it follows that $\Phi$ is a strict contraction on
$\mathcal{H}_T^2$ equipped with the norm $\|\cdot\|_\beta$ if
$\beta=4C^2+1$. Then $\Phi$ has a unique fixed point and the theorem
is proved .$\diamond$

\noindent\textbf{Proof of Theorem 2:}

Applying It$\hat{o}$'s formula from $s=t$ to $s=T$, to
$\|\bm{Y}(s)-\bm{Y}(s)'\|^2$, it follows that
\begin{eqnarray}
\|\bm{Y}(T)-\bm{Y}(T)'\|^2-\|\bm{Y}(t)-\bm{Y}(t)'\|^2&=&2\int_t^T(\bm{Y}(s-)-\bm{Y}'(s-))\cdot d(\bm{Y}(s)-\bm{Y}(s)')\nonumber\\
&&+\sum\limits_{j=1}^m\int_t^Td[\bm{Y}_j-\bm{Y}_j',\bm{Y}_j-\bm{Y}_j'](s).\nonumber
\end{eqnarray}
Taking expectations and using the relations
\begin{eqnarray}
-d(\bm{Y}(t)-\bm{Y}(t)')&=&(\bm{f}(t,\bm{Y}(t-),\bm{Z}(t))-\bm{f}'(t,\bm{Y}(t-)',\bm{Z}(t)'))dt\nonumber\\
&&+\sum\limits_{d=1}^{\infty}\sum\limits_{\bm{p}\in\mathbb{N}_d^n}
\overline{\bm{\phi}}^{\bm{p}}(t)dH^{\bm{p}}(t),\nonumber
\end{eqnarray}

\begin{eqnarray}
\begin{array}{rl}
&d[\bm{Y}_j-\bm{Y}_j',\bm{Y}_j-\bm{Y}_j'](t)\nonumber\\
=&\sum\limits_{d=1}^{\infty}\sum\limits_{\bm{p}\in\mathbb{N}_d^n}\sum\limits_{e=1}^{\infty}\sum\limits_{\bm{q}\in\mathbb{N}_e^n}
\overline{\bm{\phi}}_j^{\bm{p}}(t)\overline{\bm{\phi}}_j^{\bm{q}}(t)d[H^{\bm{p}},H^{\bm{q}}](t),\qquad
j=1,2,\cdots,m\nonumber
\end{array}
\end{eqnarray}

\begin{eqnarray}
<H^{\bm{p}},H^{\bm{q}}>(t)=\delta_{\bm{p},\bm{q}}t. \nonumber
\end{eqnarray}

we have
\begin{eqnarray}
\begin{array}{rl}
&\mathbb{E}\left[\|\bm{Y}(t)-\bm{Y}(t)'\|^2\right]
+\sum\limits_{j=1}^m\sum\limits_{d=1}^{\infty}\sum\limits_{\bm{p}\in\mathbb{N}_d^n}
\mathbb{E}\left[\int_t^T\overline{\bm{\phi}}_j^{\bm{p}}(s)^2ds\right]\nonumber\\
=&\mathbb{E}\left[\int_t^T\|\bm{\xi}-\bm{\xi}'\|^2ds\right]\nonumber\\
&+2\mathbb{E}\left[\int_t^T(\bm{Y}(s-)-\bm{Y}'(s-))\cdot(\bm{f}(s,\bm{Y}(s-),\bm{Z}(s))-\bm{f}'(s,\bm{Y}(s-)',\bm{Z}(s)'))ds\right].\nonumber
\end{array}
\end{eqnarray}
Using the Lipschitz property of $\bm{f}'$, and computations similar
to those of the proof of Theorem 1 we obtain
\begin{eqnarray}
\begin{array}{rl}
&\mathbb{E}\left[\|\bm{Y}(t)-\bm{Y}(t)'\|^2\right]
+\frac{1}{2}\mathbb{E}\left[\int_t^T\|\overline{\bm{\phi}}(s)\|^2ds\right]\nonumber\\
\leq&\mathbb{E}\left[\|\bm{\xi}-\bm{\xi}'\|^2\right]+(1+2C'+2C'^2)\mathbb{E}\left[\int_t^T\|\bm{Y}(s-)-\bm{Y}'(s-)\|^2ds\right]\nonumber\\
&+\mathbb{E}\left[\int_t^T\|\bm{f}(s,\bm{Y}(s-),\bm{Z}(s))-\bm{f}'(s,\bm{Y}(s-)',\bm{Z}(s)')\|^2ds\right].\nonumber
\end{array}
\end{eqnarray}
Then by Gronwall's inequality the result follows.$\diamond$

\begin{lem}
Let $h:\Omega\times[0,T]\times\mathbb{R}^n\rightarrow\mathbb{R}$ be
a random function measurable with respect to
$\mathcal{P}\otimes\mathscr{B}_{\mathbb{R}^n}$ such that
\begin{eqnarray}
|h(s,\bm{y})|\leq a_s(\bm{y}\cdot\bm{y}\wedge\|\bm{y}\|)\quad a.s.,
\end{eqnarray}
where $\{a_s,0\leq s\leq T\}$ is a nonnegative predictable process
such that $\mathbb{E}\left[\int_0^Ta_s^2ds\right]<\infty$. Then for
each $t\in[0,T]$ we have
\begin{eqnarray}
\sum\limits_{t<s\leq T}h(s,\triangle \bm{X}(s))
=\int_t^T\int_{\mathbb{R}^n}h(s,\bm{y})\nu(d\bm{y})ds+\sum\limits_{d=1}^{\infty}\sum\limits_{\bm{p}\in\mathbb{N}_d^n}\int_0^t
<h(s,\cdot),\textsl{p}^{\bm{p}}>dH^{\bm{p}}(s) .\nonumber
\end{eqnarray}

\end{lem}
\noindent\textbf{Proof of Lemma 5:} Because (13) implies that
$\mathbb{E}\left[\int_0^t\int_{\mathbb{R}^n}|h(s,\bm{y})|^2\nu(d\bm{y})ds\right]<\infty$,
we have that
\begin{eqnarray}
M(t)=\sum\limits_{0<s\leq t}h(s,\triangle
\bm{X}(s))-\int_0^t\int_{\mathbb{R}^n}h(s,\bm{y})\nu(d\bm{y})ds.\nonumber
\end{eqnarray}
is a square integrable martingale. By the Predictable Representation
Theorem, there exists a process $\bm{\phi}$ in the space
$(M_T^2(l^2))^{\otimes n}$ such that
\begin{eqnarray}
M(t)=&\sum\limits_{d=1}^{\infty}\sum\limits_{\bm{p}\in\mathbb{N}_d^n}\int_0^t\phi^{\bm{p}}(s)dH^{\bm{p}}(s)\nonumber
\end{eqnarray}
Taking into account that
$<H^{\bm{p}},H^{\bm{q}}>_t=t\delta_{\bm{p}\bm{q}}$, we have
\begin{eqnarray}
<M,H^{\bm{p}}>_t=\int_0^t\phi^{\bm{p}}(s)ds.
\end{eqnarray}
On the other hand, using that $\triangle M(s)\triangle
H^{\bm{p}}(s)=h(s,\triangle \bm{X}(s))\textsl{p}^{\bm{p}}(\triangle
\bm{X}(s))$ we obtain
\begin{eqnarray}
<M,H^{\bm{p}}>_t=\int_0^t\int_{\mathbb{R}^n}h(s,\bm{y})\textsl{p}^{\bm{p}}(\bm{y})\nu(d\bm{y})ds.
\end{eqnarray}
Consequently, (14) and (15) imply
\begin{eqnarray}
\phi^{\bm{p}}(s)=\int_{\mathbb{R}^n}h(s,\bm{y})\textsl{p}^{\bm{p}}(\bm{y})\nu(d\bm{y}).\nonumber
\end{eqnarray}
and the result follows.$\diamond$

\noindent\textbf{Proof of Proposition 3:}

Under the hypotheses of Proposition 3 the function
$\theta_k^{(1)}(t,\bm{x},\bm{y})$ for $k=1,2,\cdots,m$ given by (6)
satisfies the hypotheses in Lemma 5 imposed on $h$ due to the mean
value theorem, when we take $\bm{x}=\bm{X}(t-)$.

Apply It\^{o}'s lemma to $\theta_k(s,\bm{X}(s))$ from $s=t$ to
$s=T$:
\begin{eqnarray}
\theta_k(T,\bm{X}(T))-\theta_k(t,\bm{X}(t))&=&\int_t^T\frac{\partial\theta_k}{\partial
t}(s,\bm{X}(s-))ds+\sum\limits_{i=1}^n\int_t^T\frac{\partial\theta_k}{\partial
x_i}(s,\bm{X}(s-))dX_i(s)\\
&&+\sum\limits_{t<s\leq
T}\left[\theta_k(s,\bm{X}(s))-\theta_k(s,\bm{X}(s-))-\sum\limits_{i=1}^n\frac{\partial\theta_k}{\partial
x_i}(s,\bm{X}(s-))\triangle X_i(s)\right].\nonumber
\end{eqnarray}
If we apply Lemma 5 to
$h(s,\bm{X},\bm{y})=\theta_k(s,\bm{X}(s-)+\bm{y})-\theta_k(s,\bm{X}(s-))-\sum\limits_{i=1}^n\frac{\partial\theta_k}{\partial
x_i}(s,\bm{X}(s-))y_i$, we obtain
\begin{eqnarray}
\begin{array}{rl}
&\sum\limits_{t<s\leq
T}\left[\theta_k(s,\bm{X}(s))-\theta_k(s,\bm{X}(s-))-\sum\limits_{i=1}^n\frac{\partial\theta_k}{\partial
x_i}(s,\bm{X}(s-))\triangle X_i(s)\right].\\
=&\sum\limits_{d=1}^{\infty}\sum\limits_{\bm{p}\in\mathbb{N}_d^n}\int_t^T\left(\int_{\mathbb{R}^{n}}\theta_k^{(1)}(s,\bm{X}(s-),\bm{y})
\textsl{p}^{\bm{p}}(\bm{y})\nu(d\bm{y})\right)dH^{\bm{p}}(s)\\
&+\int_t^T\int_{\mathbb{R}^n}\theta_k^{(1)}(s,\bm{X}(s-),\bm{y})\nu(d\bm{y})ds
\end{array}
\end{eqnarray}
Hence, substituting (17) into (16) yields
\begin{eqnarray}
\begin{array}{rl}
&g_k(\bm{X}(T))-\theta_k(t,\bm{X}(t))\\
=&\int_t^T\frac{\partial\theta_k}{\partial
t}(s,\bm{X}(s-))ds+\sum\limits_{i=1}^n\int_t^T\frac{\partial\theta_k}{\partial
x_i}(s,\bm{X}(s-))dX_i(s)\\
&+\sum\limits_{d=1}^{\infty}\sum\limits_{\bm{p}\in\mathbb{N}_d^n}\int_t^T
\left(\int_{\mathbb{R}^{n}}\theta_k^{(1)}(s,\bm{X}(s-),\bm{y})
\textsl{p}^{\bm{p}}(\bm{y})\nu(d\bm{y})\right)dH^{\bm{p}}(s)\\
&+\int_t^T\int_{\mathbb{R}^n}\theta_k^{(1)}(s,\bm{X}(s-),\bm{y})\nu(d\bm{y})ds
\end{array}
\end{eqnarray}
Notice that
\begin{eqnarray}
X_i(t)=Y_i^{(1)}(t)+t\mathbb{E}(X_i(1))=\sum\limits_{j=1}^n\tilde{c}_{ij}H^{\bm{e}_j}(t)+t\mathbb{E}(X_i(1)),\nonumber
\end{eqnarray}
and
\begin{eqnarray}
\mathbb{E}(X_i(1))=a_i+\int_{\{|y_i|\geq
1\}}y_i\nu(d\bm{y}).\nonumber
\end{eqnarray}
By applying the condition (5) and
$\bm{Y}(0)=\mathbb{E}[\bm{Y}(0)]=\mathbb{E}[g(\bm{X}(T))]$, We can
rewrite (18) as
\begin{eqnarray}
g_k(\bm{X}(T))&=&\mathbb{E}[g_k(\bm{X}(T))]+\sum\limits_{i=1}^n\int_t^T\frac{\partial\theta_k}{\partial
x_i}(s,\bm{X}(s-))\sum\limits_{j=1}^n\tilde{c}_{ij}dH^{\bm{e}_j}(t)\nonumber\\
&&+\sum\limits_{d=1}^{\infty}\sum\limits_{\bm{p}\in\mathbb{N}_d^n}\int_t^T
\left(\int_{\mathbb{R}^{n}}\theta_k^{(1)}(s,\bm{X}(s-),\bm{y})
\textsl{p}^{\bm{p}}(\bm{y})\nu(d\bm{y})\right)dH^{\bm{p}}(s)
\end{eqnarray}
which completes the proof of the Proposition. $\diamond$

\noindent\textbf{Proof of Proposition 4:}

Apply It\^{o}'s lemma to $\theta_k(s,\bm{X}(s))$ for
$k=1,2,\cdots,d$ from $s=t$ to $s=T$. By using Lemma 5, we obtain
the equality (18). Now, using (10) we get
\begin{eqnarray}
\begin{array}{rl}
&g_k(\bm{X}(T))-\theta_k(t,\bm{X}(t))=-\int_t^Tf_k(s,\bm{Y}(s-),\bm{Z}(s))ds\\
&+\sum\limits_{i=1}^n\int_t^T\frac{\partial\theta_k}{\partial
x_i}(s,\bm{X}(s-))\sum\limits_{j=1}^n\tilde{c}_{ij}dH^{\bm{e}_j}(s)
+\sum\limits_{i=1}^n\int_t^T\left(\int_{\mathbb{R}^{n}}\theta_k^{(1)}(s,\bm{X}(s-),\bm{y})\textsl{p}^{(1)}(\bm{y})\nu(d\bm{y})\right)dH^{\bm{e}_i}(s)\\
&+\sum\limits_{d=2}^{\infty}\sum\limits_{\bm{p}\in\mathbb{N}_d^n}\int_t^T\left(\int_{\mathbb{R}^{n}}\theta_k^{(1)}(s,\bm{X}(s-),\bm{y})\textsl{p}^{\bm{p}}(\bm{y})
\nu(d\bm{y})\right)dH^{\bm{p}}(s)\nonumber
\end{array}
\end{eqnarray}
which completes the proof of the Proposition.

\noindent\textbf{Acknowledgement}: The authors would like to thank a
kind proposal given by Professor David Nualart for the initial
version of this paper.

\end{document}